\documentclass{amsart}
\usepackage{graphicx}
\newtheorem{theorem}{Theorem}

\newtheorem{conjecture}[theorem]{Conjecture}

\newcommand{\abs}[1]{\left\vert#1\right\vert}

\usepackage{amsmath,amscd}
\begin{document}
\title[]{On the Bishop Invariants Of Embeddings of $S^3$ into $\mathbb{C}^3$}%
\author{Ali M. Elgindi}%

\begin{abstract}
The Bishop invariant is a powerful tool in the analysis of real submanifolds of complex space that associates to every (non-degenerate) complex tangent
of the embedding a non-negative real number (or infinity). It is a biholomorphism invariant that gives information regarding the local hull of holomorphy
of the manifold near the complex tangent. In this paper, we derive a readily applicable formula for the computation of the Bishop invariant for graphical
 embeddings of 3-manifolds into $\mathbb{C}^3$. We then exhibit some examples over $S^3$ demonstrating the different possible configurations of the
 Bishop invariant along complex tangents to such embeddings. We will also generate a few more results regarding the behavior  of the Bishop invariant in
 certain situations.
  We end our paper by analyzing the different possible outcomes from the perturbation of a degenerate complex tangent.

   \end{abstract}

\maketitle

\section*{I. Introduction}
In this paper, our main focus will be on the subject of the Bishop invariant of complex tangents to embeddings of 3-manifolds into $\mathbb{C}^3$. We will
in particular focus on graphical embeddings of $S^3$. The Bishop invariant is of considerable interest in complex analysis as the value of the Bishop
invariant at a complex tangent determines the behavior of the hull of holomorphy near the complex tangent. This result was discovered by Errett Bishop,
which he published in [2]. Further work by Moser and Webster in 1983 (see [6]) determined precisely the normal form of a real surface embedding in $\mathbb{C}^2$
near an elliptic complex tangent, and in doing so determined precisely the local hull of holomorphy for such a surface near an elliptic tangent.
The following year, Kenig and Webster (in [4]) extended the previous work to determine the hull of holomorphy of a real $n$-manifold embedded in
$\mathbb{C}^n$.
\par\
For practical purposes however, it is unfortunate that the derivation of the Bishop invariant can be found to be rather coarse and abstract, and not easily applicable for
computation. In our paper, we will use the work of Bishop and the later work of Webster in [7] to derive a new formula for the Bishop invariant which can be used to
easily compute the Bishop invariant for graphical embeddings of 3-manifolds. We demonstrate through several examples the flexibility of the Bishop
invariant in this dimension, in particular we exhibit the possible behaviors of the Bishop invariant along different configurations of complex tangents,
which in the dimension three will generically arise as knots or links. We also prove that along any surface of complex tangents, the Bishop invariant must be $\frac{1}{2}$.
We then analyze the effect of perturbing a degenerate complex tangent, both in terms of the structure of the complex tangents to the perturbed embedding
and the corresponding behavior of the Bishop invariant.

\par\ \par\

\section*{II. On Bishop Invariants to 3-manifolds}
The $\emph{Bishop invariant}$ is a biholomorphism invariant associated to embeddings of real manifolds $M^n \hookrightarrow \mathbb{C}^n$. The Bishop invariant associates a real number (or infinity) to every non-degenerate complex tangent point $x$ of degree one, i.e. a point
$x \in M \subset \mathbb{C}^n$ so that $T_x (M) \cap J(T_x (M)) \leq \mathbb{C}^n$ is a complex line (satisfying a non-degeneracy condition to be given shortly).
\par\ \par\
More precisely, for such a complex tangent $x \in M$, there is an associated number $\gamma(x) \in [0, \infty]$, and we call $x$ $\emph{elliptic}$ if $\gamma(x) \in [0, \frac{1}{2})$, $\emph{hyperbolic}$ if $\gamma(x) \in (\frac{1}{2}, \infty]$, and $\emph{parabolic}$ if $\gamma(x) = \frac{1}{2}$. It is a celebrated
result of Bishop in complex analysis that if a complex tangent $x$ is elliptic then there exists a family of analytic disks in $\mathbb{C}^n$ whose boundary circles all lie in $M$ and so that the disks converge to the point $x$. Using Bishop's results, Webster and his colleagues later determined precisely the hull of holomorphy for an n-manifold embedded in n-dimensional complex Euclidean space. We refer the reader to Bishop's publication in [1], and the works of Webster-Moser ([6]), and Webster-Kenig ([4]).
\par\ \par\
In our paper, we wish to investigate if there is any topological or absolute restrictions on the behavior of the Bishop invariant along curves of complex tangents in $S^3$. From our computations, there seems to be many possible behaviours of $\gamma$ along any knot component. In particular, there are situations where $\gamma$ is constant along a curve, while there are situations where it varies very dramatically along a single curve and assumes every type of complex tangent (hyperbolic, parabolic, and elliptic) on the same connected component.
\par\
First, we must derive $\gamma$ and find the means for practical computation. As an alternative exposition on the Bishop invariant, we refer the reader to Webster in [7]. For practical purposes, from here on we restrict our considerations to 3-manifolds embedded in $\mathbb{C}^3$.
\par\ \par\
Suppose $\mathcal{S} \hookrightarrow \mathbb{C}^3$ can be given in holomorphic coordinates $(z, w, \zeta)$ by: $\mathcal{S} = \{\vec{\mathcal{R}}=(r_1,r_2,r_3) = 0\}$, where $r_i:\mathbb{C}^3 \rightarrow \mathbb{R}$ are smooth and $dr_1 \wedge dr_2 \wedge dr_3 \neq 0$ (nowhere zero on $\mathcal{S}$). Let $\partial r_1 \wedge \partial r_2 \wedge \partial r_3 = B dz \wedge dw \wedge d\zeta$,  where $\partial$ is the holomorphic Dolbeault operator. Hence $B:\mathcal{S} \rightarrow \mathbb{C}$ is a smooth map, in fact: $B(x) = det(\frac{\partial(r_1, r_2, r_3)}{\partial(z, w, \zeta)})$. By construction, we see that $B(x)=0$ if and only if $x$ is a complex tangent.
\par\ \par\
In this paper, we will focus on the situation where $\mathcal{S}= graph(f|_M)$ is a graph of a (sufficiently) smooth complex-valued function $f$ over $M=\{\rho(z,w)=0\} \subset \mathbb{C}^2$, where $\vec{\nabla}(\rho) \neq 0$ on $M$.  In this case, we may write:
\par\
$\mathcal{S} = \{(z,w,\zeta) | \rho(z,w)=0$ and $\zeta - f(z,w) = 0\}$.
\par\
With respect to our previous notation, we write: $r_3 = \rho(z,w)$, and $r_1 = Re(\zeta-f(z,w)), r_2 = Im(\zeta-f(z,w))$. As before: $\mathcal{S} = \{\vec{\mathcal{R}} = 0\}$, where $\vec{\mathcal{R}}=(r_1,r_2,r_3)$.
\par\
Note we may change our coordinates to: $(R, \overline{R}, \widehat{R})$ where $R = \zeta - f(z,w)$ and $\widehat{R} = \rho(z,w)$ without affecting  the above computation for $B$ (by properties of the complex determinant). Taking $M=S^3$ ($\rho(z,w) = \abs{z}^2+\abs{w}^2 -1)$, we compute:
\par\ \par\
$ B =det \left( {\begin{array}{ccc}
\frac{\partial R}{\partial z} & \frac{\partial \overline{R}}{\partial z} & \frac{\partial \widehat{R}}{\partial z} \\
\frac{\partial R}{\partial w} & \frac{\partial \overline{R}}{\partial w} & \frac{\partial \widehat{R}}{\partial w} \\
\frac{\partial R}{\partial \zeta} & \frac{\partial \overline{R}}{\partial \zeta} & \frac{\partial \widehat{R}}{\partial \zeta} \end{array}} \right) $

$= det \left( {\begin{array}{ccc}
-f_z & -\overline{f}_z & \overline{z} \\
-f_w & -\overline{f}_w & \overline{w} \\
1 & 0 & 0 \end{array} } \right) $
\par\ \par\
$=\overline{z} (\overline{f})_w - \overline{w} (\overline{f})_z = -\overline{\mathcal{L} (f)}$,
where $\mathcal{L} = z \frac{\partial}{\partial \overline{w}} - w \frac{\partial}{\partial \overline{z}}$ is the tangential Cauchy-Riemann operator to $S^3$.
\par\
Note this is consistent as then $\{B=0\}=\{\mathcal{L} (f) = 0\}$, and the zeros of $\mathcal{L} (f)$
form precisely the set of complex tangents to $\mathcal{S}= graph(f|_{S^3})$; see Ahern and Rudin in [1].
\par\ \par\
We will now follow the work of Webster in [7] in defining the Bishop invariant. Denote by $\aleph_f = \{B(x)=0\}$ the set of complex tangents to the embedding $\mathcal{S}= graph(f|_{S^3})$. Suppose $f$ is a generic map, in the sense that the set of complex tangents is a smooth manifold of codimension two in $S^3$, and so $\aleph_f$ is some disjoint union of closed curves in $S^3$, i.e. a link of knots. Then for $m \in \aleph_f$, $H_m(\mathcal{S}) = T_m (\mathcal{S}) \cap J(T_m (\mathcal{S})) \subset \mathbb{C}^3$ is a complex line. In fact, we get a complex line bundle $H \rightarrow \aleph_f$ over the link $\aleph_f$. Now, tensoring each holomorphic tangent space with $\mathbb{C}$,
we get the natural decomposition: $H_m \bigotimes \mathbb{C} = {H_m}^a \bigoplus {H_m}^b$, where ${H_m}^b = \overline{{H_m}^a}$ and ${H_m}^a$ is spanned by a vector $X = \sum \xi_i \frac{\partial}{\partial z_i}$, which is given by the condition that $X(\vec{\mathcal{R}}) = 0$ on $\aleph_f$.
\par\ \par\
We call a complex tangent $m \in M$ $\emph{non-degenerate}$ if $XB (m) \neq 0$   $\emph{\textbf{OR}}$  $\overline{X}
B (m) \neq 0$.
\par\ \par\
For such a non-degenerate complex tangent $m$,  and define the Bishop invariant at $m$ to be:
\par\
$\gamma(m) = \frac{1}{2} \abs{\frac{XB(m)}{\overline{X}B(m)}} \in [0, \infty]$.
\par\ \par\
Note that as the vector field $X$ is defined up to $X \rightarrow cX$, $c$ a complex constant, it follows that the Bishop invariant is a well-defined biholomorphism invariant.
\par\ \par\
In order to compute the Bishop invariant, we need to make sense of the vector field $X$. We may compute:
$X(\vec{\mathcal{R}}) = \xi_1 {\begin{pmatrix} -f_z \\ -(\overline{f})_z \\ \overline{z} \end{pmatrix}} +
\xi_2 {\begin{pmatrix} -f_w \\ -(\overline{f})_w \\ \overline{w} \end{pmatrix}} +
\xi_3 {\begin{pmatrix} 1 \\ 0 \\ 0
\end{pmatrix}}$
\par\ \par\
$= {\begin{pmatrix} \xi_3 - \xi_1 f_z - \xi_2 f_w \\ -\xi_1 (\overline{f})_z -\xi_2 (\overline{f})_w \\
\xi_1 \overline{z} + \xi_2 \overline{w} \end{pmatrix}}$
\par\ \par\
We wish to solve the equation: $X(\vec{\mathcal{R}})(m) = \overrightarrow{0}$ for all $m \in \aleph_f$. In particular, we need to find $\xi_1, \xi_2, \xi_3$.
(for a given function $f$) so that the above vector is zero at every point of $\aleph_f$. However, as we only need $X$ to evaluate on $B$, and $B$
is a function solely in $z, w$ (by construction), we find that $\xi_3$ is not relevant for our purposes and so we only need consider the last two equations (which are devoid of mention of $\xi_3$) to find $\xi_1, \xi_2$.
\par\
$\xi_1 (\overline{f})_z + \xi_2 (\overline{f})_w = 0$
\par\
$\xi_1 \overline{z} + \xi_2 \overline{w} = 0$, for all $(z,w) \in \aleph_F$.
\par\ \par\
Now, suppose $(z,w) \in \aleph_f$ so that $z \neq 0$. Then we may write the second equation: $\xi_1 = -\frac{\overline{w}}{\overline{z}} \xi_2$. Inserting this into the first equation, we get:
\par\ \par\
$-\frac{\overline{w}}{\overline{z}} \xi_2 (\overline{f})_z + \xi_2 (\overline{f})_w = 0$.
\par\
Multiplying both sides by $\overline{z}$, we obtain:
\par\
$\xi_2(-\overline{w} (\overline{f})_z + \overline{z} (\overline{f})_w)) = 0$.
\par\
which is now a vacuous condition as the term in the parenthesis (i.e. everything excluding $\xi_2$) is equal zero by condition that  $(z,w) \in \aleph_f$. This is in fact consistent, as there is a complex line of possible candidates for $X$, in particular any nowhere-zero complex function times a suitable $X$ will also generate the complex tangent space.
\par\ \par\
In particular, we find $X = -\frac{\overline{w}}{\overline{z}} \xi_2 \frac{\partial}{\partial z} + \xi_2 \frac{\partial}{\partial w} +\alpha \frac{\partial}{\partial \zeta}$, for any nowhere-zero (in a neighborhood of $(z,w) \in \aleph_f$) complex function $\xi_2$ and some irrelevant (though can be easily computed from $\xi_2$) function $\alpha$. Taking $\xi_2 = \overline{z}$ (remember we assume $z \neq 0$), we may take $X$ to be:
\par\ \par\
$X= -\overline{w} \frac{\partial}{\partial z} + \overline{z} \frac{\partial}{\partial w}+ \alpha \frac{\partial}{\partial \zeta}$.
\par\ \par\
Note if $z=0$ this would imply that $\abs{w} =1$ (since we are in $S^3$), and the above system of equations above would imply that $\xi_2 = 0$ which means $X= -\overline{w} \frac{\partial}{\partial z} +\alpha \frac{\partial}{\partial \zeta}$ would generate the complex tangent space. Hence, our above choice for $X$ will in fact work at any complex tangent $(z,w) \in \aleph_f$.
\par\ \par\
Our computation of interest is: $X(B)$, where $B=-\overline{\mathcal{L} (f)}$, strictly a function in $z, w$. We see now why the function $\alpha$ is irrelevant to our considerations. But once can easily see that $X$, acting on functions $g$ solely of $z,w$, acts as:
$X(g) = - \overline{w} \frac{\partial g}{\partial z} + \overline{z} \frac{\partial g}{\partial w} = -\overline{\mathcal{L}} (g)$.
\par\ \par\
Hence, we see that: $X(B) = X(-\overline{\mathcal{L} (f)}) = -\overline{\mathcal{L}} (-\overline{\mathcal{L} (f)}) =\overline{\mathcal{L}}^2 (\overline{f})$, by linearity of $X$ and by the simple fact that: $\overline{\mathcal{L} (f)} = \overline{\mathcal{L}} (\overline{f})$.
\par\ \par\
Further, we see that: $\overline{X} (B) = - \mathcal{L} (-\overline{\mathcal{L} (f)}) = \mathcal{L} (\overline{\mathcal{L}(f)}) = \mathcal{L} (\overline{\mathcal{L}} (\overline{f}))$.
\par\ \par\
In fact, we are only interested in the norms of the above functions $X(B)$ and $\overline{X}(B)$, so we may take their conjugates:
\par\
$\abs{X(B))} = \abs{\overline{X(B)}} = \abs {\mathcal{L}^2 (f)}$
\par\
$\abs{\overline{X} (B)} = \abs{\overline{\mathcal{L}} (\mathcal{L}(f))}$,
\par\
both of which are clear by our above observations.
\par\ \par\
As such, we may now give the Bishop invariant a more closed and direct formulation for graphical embeddings:
\par\
\begin{theorem} Let $G:S^3 \hookrightarrow \mathbb{C}^3$  be a "graphical embedding"; i.e. $G=graph(g|_{S^3})$, for some $g:\mathbb{C}^2 \rightarrow \mathbb{C}$ smooth (at least $\mathcal{C}^2$). Then for any complex tangent $(z,w) \in \aleph_g$, we may write the Bishop invariant as:
\par\
$\gamma_g (z,w) = \frac{1}{2} \frac{\abs{\mathcal{L}^2 (g)}}{\abs{\overline{\mathcal{L}}(\mathcal{L}(g))}} (z,w) \in [0, \infty]$, where the complex tangent is degenerate if both the numerator and dominator are zero.
\end{theorem}
\par\
We remark that this is a biholomorphic invariant as it is derived directly from the standard formulation of Bishop (in [2]) and the later work of
Webster (in [7]).
\par\ \par\
Let us now make a note of the value of the Bishop invariant at a complex tangent $x \in \aleph_f$ which lies on a surface of complex tangents.
\par\
Assume that $\aleph_f = \{ \mathcal{L} (f) = 0 \} \bigcap S^3$ contains a surface, that is, a codimension 1 submanfiold of $S^3$. We will assume the
situation that $\mathcal{L} (f)$ is a real function, or a complex function times a real function. We say this surface is of $\emph{generic type}$.
 Moreover, if $\aleph_f$ contains a surface which contains $x$ but no other complex tangent component contains $x$ (i.e. $x$ is not the intersection
 of a surface complex tangent component and a curve complex component, for example), we may write $\mathcal{L} (f) (z,w) = c r(z,w) g(z,w)$, where
 $r: \mathbb{C}^2 \rightarrow \mathbb{R}$ is real and $r(x) = 0$ but $g(x) \neq 0$, and $c \in \mathbb{C}$. We call such $x$ $\emph{non-singular}$.
\par\ \par\
We then compute:
\par\
$\mathcal{L}^2 (f) (z,w) = \mathcal{L} (\mathcal{L} (f)) (z,w) = \mathcal{L} (c r(z,w) g(z,w)) =c g(z,w) \mathcal{L} (r(z,w)) +c r(z, w) \mathcal{L} (g(z,w))$.
\par\ \par\
Evaluating at $x$:
\par\
$\mathcal{L}^2 (f) (x) = c g(x) \mathcal{L} (r) (x) +c r(x) \mathcal{L} (g) (x) = c g(x) \mathcal{L} (r) (x)$, as $r(x) =0$.
\par\ \par\
Similarly, we see that:
\par\
$\overline{\mathcal{L}}(\mathcal{L}(f)) (x)= c \overline{\mathcal{L}}(rg) (x)= c r(x) \overline{\mathcal{L}}(g)(x) + c g(x) \overline{\mathcal{L}}(r) (x) = c g(x) \overline{\mathcal{L}}(r)(x)$.
\par\ \par\
Further, as $r$ is real, we see that: $\overline{\mathcal{L}(r)} = \overline{\mathcal{L}} (r)$. Hence, we compute the Bishop invariant at $x$:
\par\ \par\
$\gamma_f (x) = \frac{1}{2} \frac{\abs{\mathcal{L}^2 (f)}}{\abs{\overline{\mathcal{L}}(\mathcal{L}(f))}} (x) = \frac{1}{2} \frac{\abs{c g(x)} \abs{\mathcal{L} (r) (x)}}{\abs{c g(x)} \abs{\overline{\mathcal{L}}(r)(x)}}$.
\par\
But we know $c g(x) \neq 0$, so they cancel and by our above remark the resulting numerator and denominator functions are conjugate (and hence have the same norm):
\par\ \par\
$\gamma_f (x) = \frac{1}{2} \frac{\abs{\mathcal{L} (r) (x)}}{\abs{\overline{\mathcal{L}}(r)(x)}} = \frac{1}{2}$.
\par\ \par\
Hence, we have proven:
\par\
\begin{theorem} Let $E:S^3 \hookrightarrow \mathbb{C}^3$ be a graphical embedding. Then any (non-singular) complex tangent lying on a surface of complex tangents of generic type must be parabolic (and in particular, with non-degenerate Bishop invariant).  \end{theorem}
\par\
Surfaces may arise as complex tangents only in non-generic situations (only for a subset of embeddings of measure zero), but arise in many approachable examples; perhaps the simplest situation is when $g(z,w) = \overline{z} \overline{w}$, where it is easy to see that the embedding under the graph of $g$ will give rise to complex tangents precisely along the torus: $\mathcal{T} =\{ (z,w) | \abs{z}^2 =\abs{w}^2 = \frac{1}{2} \} \subset S^3$, all of which are parabolic of course by our above result (there are no other components of complex tangents to need to worry about singularities).
\par\ \par\
\section*{III. Some Interesting Examples}
Let us now do some examples for the generic situation, i.e. where we have (unions) of closed curves as the complex tangents, to point out the numerous configurations and topological flexibility for the behavior of the Bishop invariants in the generic situation.
\par\ \par\
$\textbf{\emph{\underline{Ex1}}}:-$   Here we will exhibit an embedding (of $S^3$) that admits a curve of complex tangents which is purely parabolic. Note this is in general a very rare situation as it requires the Bishop invariant to take value precisely $\frac{1}{2}$ at every point of the curve.
\par\
Let $f(z,w) = \frac{1}{2} \overline{z}^2 + \overline{w}$. Then the embedding $F=graph(f|_{S^3})$ has its complex tangents precisely where $\mathcal{L} (f) = 0$.
\par\ \par\
Here $\mathcal{L} (f) = w \overline{z} -z$, and we further compute:
\par\
$\mathcal{L} (\mathcal{L} (f)) = w^2$ and $\overline{\mathcal{L}} (\mathcal{L} (f)) = -\overline{w}-\overline{z}^2$.
\par\ \par\
Now, $\aleph_f = \{\mathcal{L} (f) = 0 \} \bigcap S^3$. But note that: $\mathcal{L} (f) (z,w)= 0$ implies that $\overline{z} w = z$ , and taking the
norm of both sides, we find:
\par\
$\abs{z} \abs{w} = \abs{z}$.
\par\
There are now two cases: $\abs{z} =0$ or $\abs{z} \neq 0$. If $\abs{z} = 0$, this implies that $z=0$ and $\abs{w}=1$ on $S^3$, resulting in a great circle.
Assuming $\abs{z} \neq 0$, we may divide both sides of the above equation by $\abs{z}$ to obtain that: $\abs{w} =1$ once again. Hence, we find that the
set of complex tangents to this embedding will form exactly the great circle $\abs{w}=1$ in $S^3$. In our notation, we write:
\par\
$\aleph_f = \{(0,w) | \abs{w} =1\}$.
\par\ \par\
We compute the Bishop invariant along this curve:
\par\ \par\
$\gamma_f (0,w) = \frac{1}{2} \frac{\abs{\mathcal{L}^2 (f)}}{\abs{\overline{\mathcal{L}}(\mathcal{L}(f))}} (0,w)=
\frac{1}{2} \frac{\abs{w^2}}{\abs{-\overline{w}-\overline{z}^2}}
= \frac{1}{2} \frac{\abs{w}^2}{\abs{\overline{w}}} = \frac{1}{2}$, since $\abs{w} =1$
\par\
Hence, we see that every point on the curve $\aleph_f$ is parabolic. Hence, $\aleph_f \subset S^3$ is a great circle which is purely parabolic with respect to this particular embedding.
\par\ \par\
$\textbf{\emph{\underline{Ex2}}}:-$  We will now generalize the above example 1 and find an embedding of $S^3$ into $\mathbb{C}^3$ that admits a curve
of complex tangents with Bishop invariant constant along the curve of any prescribed value.
\par\
Let $f(z,w) = \alpha \overline{z}^2 + \overline{w}$, for any fixed $\alpha \geq 0$ (finite). We compute:
\par\
$\mathcal{L} (f) = 2 \alpha w \overline{z} -z$
\par\
$\mathcal{L} (\mathcal{L} (f)) = 2 \alpha w^2$
\par\
$\overline{\mathcal{L}} (\mathcal{L} (f)) = -\overline{w}-2 \alpha \overline{z}^2$.
\par\ \par\
First, to find the set of complex tangents: $\aleph_f = \{\mathcal{L} (f) = 0 \} \bigcap S^3$, we solve the equation:
\par\
$2\alpha w \overline{z} = z$. First, if $\alpha = 0$, we get only one curve of complex tangents corresponding to the equation $z=0$ on $S^3$.
\par\
Now, assume $\alpha > 0$. Proceeding as before, we have the two cases: $z=0$ or $z \neq 0$.
\par\ \par\
If $z=0$, we obtain the great circle: $\{(0,w) | \abs{w}=1\}$. Call this curve of complex tangents $\aleph_1$.
\par\
If $z \neq 0$, we may take the norm of both sides of the above equation and divide both
sides by $\abs{z}$ to find that: $\abs{w} = \frac{1}{2\alpha}$. This gives rise to a second curve of complex tangents, which we denote by: $\aleph_2$.
\par\ \par\
Let us consider the first curve, that is $\aleph_1$. To compute the Bishop invariant along this curve:
\par\
$\gamma_f (0,w) = \frac{1}{2} \frac{\abs{\mathcal{L}^2 (f)}}{\abs{\overline{\mathcal{L}}(\mathcal{L}(f))}} (0,w)=
\frac{1}{2} \frac{\abs{2\alpha w^2}}{\abs{-\overline{w}-2\alpha \overline{z}^2}}
= \frac{1}{2} \frac{2\alpha \abs{w}^2}{\abs{\overline{w}}} = \alpha$, since $\abs{w} =1$ and $z=0$.
\par\ \par\
Hence, $\aleph_1$ is a curve of complex tangents to the embedding given by $f$ whose Bishop invariant is constant at $\alpha$, which was an arbitrary
positive number.
\par\ \par\
We may also compute the Bishop invariant along the curve $\aleph_2$. With some direct (but somewhat messy) computation, we find that:
\par\
$\gamma_f = \frac{1}{4\alpha}$ on $\aleph_2$.
\par\ \par\
In summary, if $0<\alpha<\frac{1}{2}$, we obtain an elliptic curve $\aleph_1$ and a hyperbolic curve $\aleph_2$ of constant Bishop invariant as described above.
\par\
If $\frac{1}{2} <\alpha< \infty$, we get a hyperbolic curve $\aleph_1$ and an elliptic curve $\aleph_2$ of constant Bishop invariant.
\par\
If $\alpha = 0$, we get a single curve of elliptic complex tangents with Bishop invariant 0.
\par\
If $\alpha = \frac{1}{2}$, we get the situation of Example 1, where we have exactly one curve of complex tangents which is purely parabolic.
\par\ \par\
Finally, to be complete, we find an embedding of $S^3$ which admits complex tangents which are hyperbolic of Bishop invariant $\infty$.
\par\
Let $f(z,w) = \overline{w}^2$. Then:
\par\
$\mathcal{L}(f) = -2z \overline{w}$, and hence the set of complex tangents to this embedding will form two great circles:
\par\
$\aleph_1 = \{(0,w) | \abs{w}=1\}$ and $\aleph_2 = \{(z,0) | \abs{z}=1\}$. We further compute:
\par\
$\mathcal{L} (\mathcal{L}(f)) = 2 z^2$ and $\overline{\mathcal{L}} (\mathcal{L}(f)) = -2 \overline{w}^2$.
\par\
The Bishop invariant to this embedding at a complex tangent is given by: $\gamma_f (z,w)= \frac{1}{2} \frac{\abs{z}^2}{\abs{w}^2}$.
\par\ \par\
Hence, $\aleph_1$ consists of complex tangents which are elliptic of constant Bishop invariant 0, and $\aleph_2$ consists of complex tangents which
are hyperbolic of constant Bishop invariant $\infty$.
\par\ \par\
$\textbf{\emph{\underline{Ex3}}}:-$    Here we exhibit a collection of five great circles that form the set of complex tangents (all non-degenerate) to some embedding of $S^3$. Let $f(z,w) = \overline{z}^2 w + \overline{w}^2 z$. We compute:
\par\ \par\
$\mathcal{L} (f) = 2 \overline{z} w^2 -2 \overline{w} z^2$, and further:
\par\
$\mathcal{L} (\mathcal{L} (f)) = 2w^3 +2z^3, \overline{\mathcal{L}} (\mathcal{L} (f)) = -4(\overline{w}^2 z + \overline{z}^2 w)$.
\par\ \par\
We note that: $\aleph_a =\{(0,w)\}, \aleph_b =\{(z,0)\} \subset \aleph_f$, and as
$\mathcal{L} (\mathcal{L} (f)) \neq 0, \overline{\mathcal{L}} (\mathcal{L} (f)) = 0$ at all points of either set, we see that these two components consist purely of complex tangents of hyperbolic type-$\infty$.
\par\ \par\
Now suppose $z, w \neq 0$. It is direct to see that: $z^3 = \frac{1-\abs{w}^2}{\abs{w}^2} w^3$. Taking norm (squared) of both sides we get a polynomial in $x= \abs{w}^2$:
\par\
$x^2(2x^3-5x^2+4x-1) = 0$, whose only solutions are again: $x =0, \frac{1}{2}, 1$. Note that the solutions $x=0, x=1$ correspond exactly to $\aleph_b, \aleph_a$, respectively.
\par\ \par\
If $\abs{w}^2 =\frac{1}{2}$, we see that we get three additional curves of complex tangents:
\par\
$\aleph_1 = \{(z,z) | \abs{z}^2 =\frac{1}{2}\}, \aleph_2 = \{(z,ze^{2\pi i/3}) | \abs{z}^2 =\frac{1}{2}\}, \aleph_3 = \{(z,ze^{4\pi i/3}) | \abs{z}^2 =\frac{1}{2}\}$.
\par\ \par\
Let $(z,z) \in \aleph_1$. Then $\gamma_f (z,z) = \frac{1}{2} \frac{\abs{2w^3 +2z^3}}{\abs{4(\overline{w}^2 z + \overline{z}^2 w)}} =
\frac{1}{4} \frac{\abs{2z^3}}{\abs{2z \abs{z}^2}} = \frac{1}{4}$.
\par\ \par\
Hence, every complex tangent in $\aleph_1$ is elliptic with Bishop invariant $\frac{1}{4}$. Further, it is easy to see that complex tangents $\aleph_2, \aleph_3$ has the same Bishop invariants as the complex units factor out of the denominator and have no effect (norm =1).
\par\ \par\

$\textbf{\emph{\underline{Ex4}}}:-$  We now consider Bishop invariants along torus knots. For more information on torus knots, we refer the reader to Milnor in [5]. Let $\mathcal{T}_{p,q} = \{z^p = w^q\} \bigcap S^3$ be the standard $(p,q)$-torus knot, where $(p,q)$ are relatively prime natural numbers. We may readily graph $\mathcal{T}_{2,3}$ using software like Mathematica (via stereographic projection).
\par\ \par\
Note that $\mathcal{T}_{p,q}$ arises precisely as the set of complex tangents to the graph of $f(z,w) = z^{p-1} \overline{w} +w^{q-1} \overline{z}$,
i.e. $\mathcal{L} (f) = z^p - w^q$. But this is a holomorphic function, so $\mathcal{L} (\mathcal{L} (f)) =0$ and further:
$\overline{\mathcal{L}} (\mathcal{L} (f)) = p \overline{w} z^{p-1} + q \overline{z} w^{q-1}$ which with a parametrization for $\mathcal{T}_{p,q}$ one can show is never zero for $(z,w) \in \mathcal{T}_{p,q}$. Hence, every point of $\mathcal{T}_{p,q}$ is an elliptic complex tangent with Bishop invariant 0.
\par\ \par\
One can extend this to general complex tangent knots arising as intersections of holomorphic algebraic sets in $\mathbb{C}^2$ intersect with $S^3$ to see that if a complex tangent is non-degenerate, it must have Bishop invariant 0.
\par\ \par\
$\textbf{\emph{\underline{Ex5}}}:-$     Here we exhibit an embedding whose complex tangents make two (knot) components who link (doubly), yet one curve consists of elliptic complex tangents, and the other consists of hyperbolic complex tangents.
\par\ \par\
Let $f(z,w) = z \overline{w}^2 +\overline{z}  \overline{w}$. We compute:
\par\ \par\
$\mathcal{L} (f) = 2 \overline{w} z^2 + 2 \abs{z}^2 -1$, and so we see that $\mathcal{L} (f) = 0$ if and only if:
\par\
$w = \frac{1-2\abs{z}^2}{2 \overline{z}^2}$. Taking norms of both sides, we get the polynomial equation in $x (= \abs{z}^2)$:
\par\
$4x^3 -4x+1 = 0$, whose only (relevant) solutions are: $x =.27, .838$. Let $\alpha_1 = .27, \alpha_2 = .838$. Then it is clear that:
\par\ \par\
$\aleph_f = \aleph_1 \bigcup \aleph_2$, where $\aleph_i = \{(z, \frac{1-2\alpha_i}{2 \overline{z}^2}) | \abs{z}^2 = \alpha_i\}$.
\par\ \par\
Using stereographic projection through the point $(0,1) \in S^3$ (which both curves avoid), we may plot the curves $\aleph_1, \aleph_2$  in $\mathbb{R}^3$, and there-in it is clear that $\aleph_1, \aleph_2$ are topologically equivalent to circles and link together with linking number 2 (using Mathematica).
\par\
In particular, this gives us an example of a non-simple link arising precisely as the set of complex tangents to an embedding.
\par\ \par\
Now, let's compute the Bishop invariants of the complex tangents. We will investigate both curves simultaneously by considering $\alpha$ in general.
\par\ \par\
Over $\aleph_f$, we may reduce: $\mathcal{L} (f) = 2 \alpha -1 +2 \overline{w} z^2$, and so:
\par\
$\mathcal{L} (\mathcal{L} (f)) = -2z^3, \overline{\mathcal{L}} (\mathcal{L} (f)) = 4 \abs{w}^2 z = 4(1- \alpha) z$.
\par\ \par\
Hence, $\gamma_f (z,w) = \frac{1}{2} \frac{\abs{-2z^3}}{\abs{4(1- \alpha) z}}= \frac{\alpha}{4(1-\alpha)}$.
\par\ \par\
On $\aleph_1, \alpha \approx .27$ and hence $\gamma_f \approx .09$ on $\aleph_1$, i.e. $\aleph_1$ consists of elliptic complex tangents.
\par\
On $\aleph_2, \alpha \approx .838$ and hence $\gamma_f \approx 1.29$ on $\aleph_2$, i.e. $\aleph_2$ consists of hyperbolic complex tangents.
\par\ \par\
The flexibility demonstrated through the examples above should hardly be a surprise, as the Bishop invariant is a local biholomorphism invariant and we expect the behavior of $\gamma$ on one component to be independent of any other component.
\par\
In fact, as we can construct unions of components of complex tangents via multiplication of relevant polynomials (at least in $\mathbb{H}$, which we may then "push-forward" to $S^3$ per our previous constructions), and if we take $\mathcal{L} (f) = g h$, where $g(x) = 0, h(x) \neq 0$, then in the computation of $\gamma_f (x)$ any mention of $h$ cancels out:
\par\ \par\
$\mathcal{L}^2 (f) (z,w) = \mathcal{L} (\mathcal{L} (f)) (z,w) = \mathcal{L} (h(z,w) g(z,w)) = g(z,w) \mathcal{L} (h(z,w)) +h(z, w) \mathcal{L} (g(z,w))$.
\par\ \par\
Evaluating at $x$:
\par\ \par\
$\mathcal{L}^2 (f) (x) = g(x) \mathcal{L} (h) (x) +h(x) \mathcal{L} (g) (x) = g(x) \mathcal{L} (h) (x)$, as $h(x) =0$.
\par\ \par\
Similarly, we see that:
\par\
$\overline{\mathcal{L}}(\mathcal{L}(f)) (x)= \overline{\mathcal{L}}(gh) (x)= h(x) \overline{\mathcal{L}}(g)(x) + g(x) \overline{\mathcal{L}}(h) (x) = g(x) \overline{\mathcal{L}}(h)(x)$.
\par\ \par\
Hence, we compute the Bishop invariant at $x$:
\par\ \par\
$\gamma_f (x) = \frac{1}{2} \frac{\abs{\mathcal{L}^2 (f)}}{\abs{\overline{\mathcal{L}}(\mathcal{L}(f))}} (x) = \frac{1}{2} \frac{\abs{g(x)} \abs{\mathcal{L} (h) (x)}}{\abs{g(x)} \abs{\overline{\mathcal{L}}(h)(x)}} = \frac{1}{2} \frac{\abs{\mathcal{L} (h) (x)}}{\abs{\overline{\mathcal{L}}(h)(x)}}$,
\par\
since we know that $g(x) \neq 0$ and as such the terms $\abs{g(x)}$ cancel out. Hence, the function $g$ has no effect on the computation of the Bishop invariant. As a result, we may conclude that the Bishop invariant is only dependent on the function that gives rise to the component of complex tangents that contains x. Hence, we have proven (at least for graphical embeddings)

\begin{theorem} Let $F=graph(f)$ be a graphical embedding. Then behavior of the Bishop invariant on any component of complex tangents is independent of the behavior on other components and only dependent on the summand of $\mathcal{L} (f)$ that gives rise to the particular component. \end{theorem}
In light of the theorem and our examples above, we may remark (a bit imprecisely):
\par\ \par\
$\emph{\textbf{\underline{Remark}}}$ The linking of two curve components of complex tangents puts no restriction on the behavior of the Bishop invariant along the two curves.
\par\ \par\
$\textbf{\emph{\underline{Ex6}}}:-$    We have thus far given examples of embeddings who give rise to knots as complex tangents, links of knots, surfaces, and unions there-in (even singular). Here we give an example of an embedding of $S^3$ who has precisely one complex tangent point.
\par\ \par\
Let $f(z,w) = \overline{z}$. We note that $H=graph(f|_\mathbb{H})$ gives a totally real embedding of $\mathbb{H}$. Consider $f \circ \psi : S^3 \setminus \{(0,1)\} \rightarrow \mathbb{C}$, in particular: $f(\psi(z,w))= \frac{i\overline{z}}{\overline{w}-1}$.
\par\ \par\
Let $q = {q_f}^3 =  \frac{i\overline{z}(1-w)^3}{\overline{w}-1}$ which is now continuously differentiable on $S^3$ (in particular, at $(0,1)$).
\par\ \par\
We are assured by our work in [3] that $Q = graph (q|_{S^3})$ will give an embedding of $S^3$ of class $\mathcal{C}^1$ with sole complex tangent at $(0,1)$. For the sake of completeness, we check:
\par\
$\mathcal{L} (q) = -iw \frac{(1-w)^3}{1-\overline{w}} + i \abs{z}^2 \frac{(1-w)^3}{(1-\overline{w})^2}$. Note both terms are continuous at $(0,1)$ and take value 0 there, in particular $(0,1)$ is a complex tangent of the embedding. Now to check there are no more complex tangents, assume we are away from $w=1$.
\par\ \par\
Then $\mathcal{L} (q) = 0$ if and only if: $\abs{z}^2 -w(1-\overline{w}) = 0$ if and only if: $w= \abs{z}^2+\abs{w}^2 = 1$, which contradicts our above assumption that we are away from $(0,1)$.
\par\ \par\
Hence, the only complex tangent to the $\mathcal{C}^1$ embedding $Q$ is $(0,1)$. As such, there exists embeddings of $S^3 \hookrightarrow \mathbb{C}^3$ who have complex tangents only at a single point.
\par\ \par\
Further, a direct computation gives us that the Bishop invariant at this point is degenerate, i.e. the dominator and numerator must both be zero; see below.
\par\
We were unsuccessful in proving the degeneracy of the Bishop invariant at a general isolated complex tangent to an embedding, however our attempts left us with confidence with which we claim:
\begin{conjecture} An isolated complex tangent to an embedding of a 3-manifold into $\mathbb{C}^3$ must have degenerate Bishop invariant.       \end{conjecture}
\par\ \par\
We now refer the reader to our paper in [3] where we derived certain functions to construct (graphical) embeddings $S^3 \hookrightarrow \mathbb{C}^3$ taking complex tangents along a knot (link) of any prescribed type. As we found, all our constructed embeddings must assume a degenerate complex tangent at the point $(0,1)$. Namely, consider maps:
\par\
$f(z,w) = \frac{(1-w)^{n+r}}{(1-\overline{w})^n} p(z,w)$  for some complex polynomial $p$ of degree $n$. Then we computed:
\par\
$\mathcal{L}({f}) (z,w) = \mathcal{L}(\frac{(1-w)^{n+r}}{(1-\overline{w})^n} p(z,w)) = \frac{(1-w)^{n+r}}{(1-\overline{w})^n} \mathcal{L}(p) (z,w) - n \frac{(1-w)^{n+r}}{(1-\overline{w})^{n+1}} z p(z,w) $, which is zero at the point $(0,1)$.
\par\ \par\
In the computations of $\mathcal{L} (\mathcal{L} (f)), \overline{\mathcal{L}} (\mathcal{L} (f))$, we find that:
\par\
$\mathcal{L} (\mathcal{L} (f)) = \frac{(1-w)^{n+r}}{(1-\overline{w})^{n+2}} g(z,w)$ for some polynomial $g$, and similarly for $\overline{\mathcal{L}} (\mathcal{L} (f))$. Taking $r >2$, we see that both must evaluate to be zero at the point $(0,1)$. Note we must take $r >2$ for the Bishop invariant to be defined.
\par\ \par\
Hence our derived embeddings necessarily come with degeneracy at the point $(0,1)$. However, by Sard's theorem, perturbing our maps $f(z,w)$ will generically (almost always) give rise to an embedding with solely non-degenerate complex tangents, however in doing so we may add to the topology of the tangents. We may also keep the perturbation "small", so that we may assume (again by Sard's theorem) that the structure of the complex tangents away from an arbitrarily small neighborhood of $(0,1)$ is not affected.
\par\ \par\
Let us now do some examples where we perturb functions of the prescribed type above and see what happens to the set of complex tangents. Note that holomorphic terms are ignored by the CR-operator so we restrict our attention to polynomials of purely non-holomorphic terms. In particular, let us continue the above example 6.
\par\ \par\
$\textbf{\emph{\underline{Ex6b}}}:-$ Here for simplicity write: $f(z,w) = \overline{z} \frac{(1-w)^4}{1-\overline{w}}$. Then as we showed above, $graph(f)|_{S^3} \subset \mathbb{C}^3$ assumes precisely one complex tangent point at $(0,1)$. Take the perturbation:
\par\
$\widetilde{f}(z,w) = \overline{z} \frac{(1+ \epsilon - w)^4}{1+\epsilon - \overline{w}}$, for $\epsilon >0$ arbitrarily small.
\par\ \par\
Then we compute (and simplify):
\par\
$\mathcal{L}(\widetilde{f})(z,w) = -\abs{z}^2  \frac{(1+ \epsilon - w)^4}{(1+\epsilon - \overline{w})^2} - w \frac{(1+ \epsilon - w)^4}{1+\epsilon - \overline{w}}  = \frac{(1+ \epsilon - w)^4}{(1+\epsilon - \overline{w})^2}  (-\abs{z}^2 -w(1+\epsilon - \overline{w}) =
\frac{(1+ \epsilon - w)^4}{(1+\epsilon - \overline{w})^2} (-\abs{z}^2 -\abs{w}^2 + (1+\epsilon)w) =
\frac{(1+ \epsilon - w)^4}{(1+\epsilon - \overline{w})^2} ((1+\epsilon)w -1)$.
\par\ \par\
Note that the fraction term is only zero for $w=1+\epsilon$, which does not lie on the sphere $S^3$. However, the other term is zero precisely when $w = \frac{1}{1+\epsilon}$ which does lie on the sphere, in fact their intersection will give a circle.
\par\
Hence, perturbing the function $f$ will give us an embedding whose complex tangents occur along a circle arbitrarily close to the point $(0,1)$.
\par\ \par\
Now, to check that this circle consists of non-degenerate complex tangents:
\par\ \par\
$\mathcal{L} (\mathcal{L} (f)) = 2z ((1+\epsilon)w -1) \frac{(1+ \epsilon - w)^4}{(1+\epsilon - \overline{w})^3} = 0$ when $w = \frac{1}{1+ \epsilon}$
\par\ \par\
$\overline{\mathcal{L}} (\mathcal{L} (f)) = \overline{z}( -4((1+\epsilon)w -1) \frac{(1+ \epsilon - w)^3}{(1+\epsilon - \overline{w})^2} + (1+\epsilon) \frac{(1+ \epsilon - w)^4}{(1+\epsilon - \overline{w})^2})$.
\par\
Now, the first term is once again zero for $w = \frac{1}{1+ \epsilon}$. The second term however, is never zero for such points $(z,\frac{1}{1+\epsilon})$ on the unit sphere. Hence, $\overline{\mathcal{L}} (\mathcal{L} (f)) \neq 0$ for every complex tangent point, and hence we have perturbed f to obtain an embedding with complex tangents along a "small" circle near $(0,1)$ all of which are non-degenerate.
\par\ \par\
$\textbf{\emph{\underline{Ex7}}}:-$ Here we exhibit a perturbation that leaves all the complex tangents the same, yet removes the degeneracy from the point $(0,1)$. Consider $f(z,w) = \overline{w} \frac{(1-w)^4}{1-\overline{w}}$.
\par\
Then: $\mathcal{L}(f)(z,w) = z (\frac{(1-w)^4}{(1-\overline{w})^2} \overline{w} + \frac{(1-w)^4}{1-\overline{w}}) =z \frac{(1-w)^4}{(1-\overline{w})^2} (\overline{w} +1 - \overline{w}) = z \frac{(1-w)^4}{(1-\overline{w})^2}$
\par\ \par\
Hence, the complex tangents will be where $z=0$ and $w=1$, but on the unit sphere $S^3$ this just translates to the great circle $\{(0,w)| \abs{w} = 1\}$. The point $(0,1)$ is necessarily degenerate. Now, let's perturb $f$ to get a new function:
\par\
$\widetilde{f}(z,w) = \overline{w} \frac{(1+\epsilon-w)^4}{1+\epsilon-\overline{w}}$ for any sufficiently small $\epsilon >0$. It is easy to see that:
\par\
$\mathcal{L} (\widetilde{f}) (z,w) = z (1+ \epsilon)\frac{(1+\epsilon-w)^4}{(1+\epsilon-\overline{w})^2}$.
\par\ \par\
The fractional expression is never zero on $S^3$, but we get the set of complex tangents to be the points where $z=0$, the same set as we found for (the unperturbed) $f$. Hence, a small perturbation of $f$ did not change the set of complex tangents.
\par\ \par\
We compute:
\par\
$\overline{\mathcal{L}} (\mathcal{L} (f)) = 4 \overline{z} (1+\epsilon) z \frac{(1+\epsilon-w)^3}{(1+\epsilon-\overline{w})^2}- \overline{w} (1+\epsilon) \frac{(1+\epsilon-w)^4}{(1+\epsilon-\overline{w})^2}$.
\par\ \par\
Note that the first term is zero for all points $\{(0,w)\} \in S^3$, but the second term is never zero for any such points as $\abs{w}=1$.
\par\
Hence, all the complex tangents of $graph(\widetilde{f})|_{S^3})$ are non-degenerate, and the set of complex tangents was unchanged by our perturbation. So effectively, we were able to "remove" the degeneracy in this situation without adding any new complex tangents.
\par\ \par\
$\textbf{\emph{\underline{Ex8}}}:-$ Let $f(z,w) = \frac{(1-w)^5}{(1-\overline{w})^2} \overline{z} \overline{w}$. We compute:
\par\ \par\
$\mathcal{L}(f)(z,w) = \abs{z}^2 (\frac{(1-w)^5}{(1-\overline{w})^2} + \frac{2 \overline{w} (1-w)^5}{(1-\overline{w})^3}) - \abs{w}^2 \frac{(1-w)^5}{(1-\overline{w})^2}$
\par\
$=\frac{(1-w)^5}{(1-\overline{w})^3} (\abs{z}^2 (2 \overline{w} + (1-\overline{w})) - \abs{w}^2 (1-\overline{w})) = \frac{(1-w)^5}{(1-\overline{w})^3} (\overline{w} +\abs{z}^2 -\abs{w}^2)$
\par\
$ = \frac{(1-w)^5}{(1-\overline{w})^3} (\overline{w} - 2\abs{w}^2 +1)$
\par\ \par\
The fractional term is zero at the point $(0,1)$. Now assume the other term is zero. Then:
\par\
$\overline{w} = 2\abs{w}^2-1$, which implies $w$ is real. So:
\par\
$2w^2-w-1 =0$ which implies that $w=1, -\frac{1}{2}$. $w=1$ again refers to the point $(0,1)$ on the unit sphere and the condition $w =-\frac{1}{2}$ gives a circle. Hence, the complex tangents of $graph(f)|_{S^3}$ consist of the (degenerate) point $(0,1)$ and a circle away from the point.
\par\ \par\
Now, let's perturb: $\widetilde{f}(z,w) = \overline{z} \overline{w} \frac{(1+ \epsilon-w)^5}{(1+\epsilon-\overline{w})^2}$
\par\
We may similarly compute:
\par\
$\mathcal{L}(\widetilde{f})(z,w) = \abs{z}^2 (\frac{(1+\epsilon-w)^5}{(1+\epsilon-\overline{w})^2} + \frac{2 \overline{w} (1+\epsilon-w)^5}{(1+\epsilon-\overline{w})^3}) - \abs{w}^2 \frac{(1+\epsilon-w)^5}{(1+\epsilon-\overline{w})^2}$
\par\
$= \frac{(1+\epsilon-w)^5}{(1+\epsilon-\overline{w})^3} (\abs{z}^2(\overline{w} +1 +\epsilon) -\abs{w}^2 (1+\epsilon - \overline{w})) =\frac{(1+ \epsilon-w)^5}{(1+\epsilon-\overline{w})^3} (\overline{w} + (1+ \epsilon) (1-2\abs{w}^2))$
\par\
Again, the fractional expression is never zero on the unit sphere, so let's set the other term equal zero. Then:
\par\
$\overline{w} = (1+\epsilon)(2\abs{w}^2-1)$, so $w$ must again be real. We get the quadratic equation in $w$:
\par\ \par\
$2w^2 - \frac{1}{1+\epsilon} w -1 =0$. By the quadratic formula in $w$, the solutions are given by:
\par\
$w=\frac{\frac{1}{1+\epsilon} \pm \sqrt{\frac{1}{(1+\epsilon)^2} +8}}{4} = \frac{1 \pm \sqrt{8(1+\epsilon)^2 +1}}{4(1+\epsilon)}$
\par\ \par\
For $ \epsilon >0$ small, this will give two solutions for $w$, one less than (but close to 1) and one near the value of -.5. For example, if we take $\epsilon = .1$, then the two values are given by: $w= .97, -.52$. Hence, the set of complex tangents to such a perturbation will consist of two circles $ w= \alpha, w =\beta$ on $S^3$. These circles cannot possible be linking as even their interiors share no common points.
\par\ \par\
Now we claim that both these circles are completely non-degenerate. We compute:
\par\
$\mathcal{L} (\mathcal{L} (f)) (z,w) = z (\frac{3(1+\epsilon-w)^5}{(1+\epsilon-\overline{w})^4}(\overline{w} + (1+\epsilon)(1-2\abs{w}^2)) + \frac{(1+\epsilon-w)^5}{(1+\epsilon-\overline{w})^3} (1-2w(1+\epsilon)))$
\par\ \par\
The first term is equal zero at any complex tangent, but as for the second term we see that as the fractional term is never zero, the only possible zeros would have to satisfy: $w = \frac{1}{2(1+\epsilon)}$. However, such a point could never be complex tangent since if we insert this value for w into the relevant quadratic equation (given above), we find that:
\par\
$\frac{1}{2(1+\epsilon)} + (1+\epsilon) (1 - \frac{1}{2(1+\epsilon)^2}) = 0$ which would imply that $1 + \epsilon =0$ which could never hold for any reasonable choice of $\epsilon$. Hence, $\mathcal{L} (\mathcal{L} (f)) (z,w) \neq 0$ at every complex tangent, and as such all complex tangents are non-degenerate.
\par\ \par\
For more complicated examples it is often difficult to even find the set of complex tangents to the unperturbed embedding. Furthermore it is not a guarantee that any arbitrary perturbation would lead to solely non-degenerate complex tangents.
\par\ \par\
The removal of this degeneracy without adding new complex tangents does not seem to be possible in every case, although as we exhibited above in some situations it can be done via certain perturbations.
\par\ \par\

\end{document}